\documentclass[12pt]{article}
\usepackage[centertags]{amsmath}
\usepackage{graphicx}
\usepackage{amssymb}
\usepackage{indentfirst}
\newtheorem{thm}{Theorem}
\newtheorem{prop}[thm]{Proposition}

\newtheorem{lem}{Lemma}[section]
\newtheorem{rem}[thm]{Example}

\def\C{{\cal C}}

 \def\G{{\cal G}}

\def\({\left(} \def\){\right)}
\def\srw{the simple random walk}
\begin{document}
\title {On the Speed of Random Walks on a
 Percolation Cluster of Trees\footnote{Supported in part by Grant G1999075106 from
the Ministry of Science and Technology of China.} }
\date{}

\author{Dayue CHEN \and Fuxi ZHANG}

\maketitle

 \begin{abstract}
We consider the simple random walk on the infinite cluster of the
Bernoulli bond percolation of trees, and investigate the relation
between the speed of the simple random walk and the retaining
probability by studying three classes of trees. A sufficient
condition is established for  Galton-Watson trees.
\end{abstract}

\section{Introduction}

The simple random walk $\{X_n\}$ on  graph $\G = (V, E)$ is
defined as a Markov chain on the set $V$ of vertices with
transition probability $p(x, y) = 1/ d_x$ if $x, y \in V$ and
there is an edge between them. The degree $d_x$ of vertex $x$ is
the number of edges connecting $x$ to other vertices.

Let $o$ be a fixed vertex and $|x|$ the graphic distance between
$o$ and $x$, i.e., the minimum number of edges in a path from $o$
to $x$.  Suppose  the simple random walk starts from $o$, i.e.,
$X_0 = o$. We call $\lim_n  |X_n|/n$, if it exists, the speed of
the simple random walk.  For the simple random walk in $Z^d$, the
speed is zero.  On the other hand the speed is $(d-1)/(d+1)$ for
the simple random walk on the regular tree $T_d$. By regular tree
$T_d$, we mean a tree that the degree of every vertex except the
root  is $d+1$.  The degree of the root  $o$ is $d$.

A non-trivial example is \srw\ on a Galton-Watson tree.  Let
$\{p_n, n= 0, 1, 2, \cdots\}$ be the offspring distribution of  a
Galton-Watson process.  Each individual produces offsprings
independently according to the same law.   There is a rooted tree
for each realization, and the correspondence induces a probability
measure in the set of rooted trees. A rooted tree, drawn according
to this measure, is called a Galton-Watson tree. It is shown in
\cite{LPP} that the speed on a Galton-Watson tree is a constant
a.s..  More remarkably, an explicit formula is given when $p_0 =
0$.
 \begin{equation} \label{Eq3}
\mbox{Speed } = \sum_{k=1}^\infty p_k \frac {k-1} {k+1}.
\end{equation}

In $p$-{\it Bernoulli bond percolation} in $\G$, each edge of $\G$
is independently declared {\it open} with probability $p$ and {\it
closed} with probability $1-p$. Thus a bond percolation $\omega$
is a random subset of $E$. We usually identify the percolation
$\omega$ with the subgraph of $\G$ consisting of all open edges
and their end-vertices. A connected component of this subgraph is
called an {\it open cluster}, or simply a {\it cluster}. Let
$p_c=p_c(\G)=\inf\{ p$:  there is an infinite cluster $a.s. \}$.
By coupling, the infinite cluster is increasing in $p$. The
probability that there is an infinite cluster is monotone in $p$.
When $p \in (p_c,1)$, with positive probability the open cluster
$\C$ that contains $o$ is infinite.

The simple random walk on the infinite cluster of $Z^d$ is first
studied in \cite{GKZ}. A more systematic investigation is
initiated in \cite{BLS}. We like to investigate the relation
between the speed of the simple random walk on an infinite cluster
and the retaining probability $p$ by observing three examples: the
regular tree, Galton-Watson trees, and the binary tree with pipes,
which
is obtained  by adding a pipe to each vertex of the binary tree.
See Figure 1.

\smallskip
\begin{prop} The speed of the simple random walk on an infinite
cluster of a regular tree is increasing in $p$.
\end{prop}

The main idea of the proof is to use (1) and a decomposition of an
infinite Galton-Watson tree as a {\it backbone} and {\it bushes}.
As $p$ increases, the backbone gets larger, the mean size and the
number of bushes get smaller. It is tempting to think that the
monotone relation holds for a large class of graphs.

\noindent {\bf Question: \it Is the speed of the simple random
walk on an infinite cluster of a transitive graph increasing in
$p$?}

We believe the answer is \lq\lq yes". However, this statement can
not be pushed any further, as we see in the following example of
the binary tree with pipes.  The binary tree is chosen because the
extinction probability can be explicitly calculated.  Assume that
$p > 1/2$ and the cluster containing the root $o$ is infinite. In
addition to the backbone and bushes, there are also pipes of
random length. The distribution of
the length 
is geometric and the mean length of a pipe is $p/(1-p)$.
Therefore, the larger the $p$ is, the longer a pipe is, and the
longer the time for the excursion on a pipe is.
As $p$ increases, there is a competition between the increase of
time spent on the pipe and the decrease of time spent on the bush.

\smallskip
\begin{rem}  The speed of the simple random walk on an infinite
cluster of the binary tree with pipes is
 $$
 \frac 1 3 \ \frac{(2p-1)^2  }{ p^2 +
(1-p)^2}\ \ \frac{1-p }{ \big(2 p^3 -6 p^2 + 3 p + 3 \big) },
 $$
which is not monotone in $p$.
\end{rem}
\begin{figure}[hp]
\centering

\includegraphics[scale=0.2]{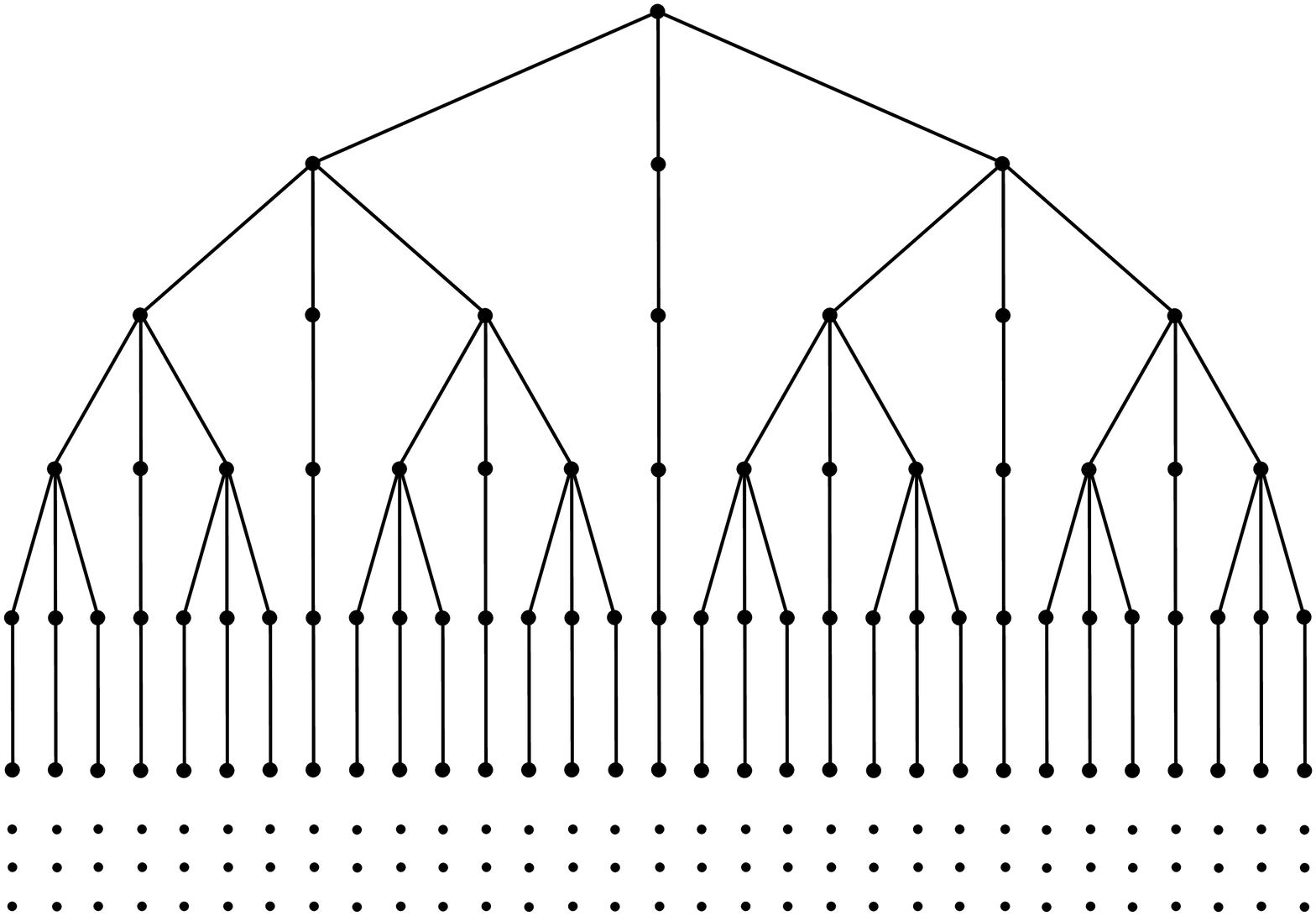}

\small{Figure 1: The binary with pipes }
\end{figure}

 Very often a random
environment (or graph) exhibits nice characters of a homogeneous
environment (or a transitive graph). So we take the Galton-Watson
trees. Let $\{p_k\}$ be the offspring distribution, and $m =
\sum_k k p_k$ the mean number of offsprings.  We assume that $m
>1$. Suppose the retaining probability $p> 1/m$ and the open
cluster $\C_o$ of $p$-Bernoulli bond percolation containing the
root is infinite. Run \srw\ starting from the root.
\begin{thm}  Let $f(s) = \sum_k p_k s^k$  be the generating function and $f'{(s)}$
the derivative of $f$. Suppose that
\begin{equation}\label{condition}
\frac {(1-s) {f'(s)}}{1-f(s) }  \mbox{ is  increasing in }\ s
\end{equation}
for $s\in (1/m, 1)$. Then the speed of the simple random walk on
an infinite cluster of a Galton-Watson tree is an increasing
function of $p$. Furthermore it is continuous for $p\in [1/m, 1]$
and differentiable for $p\in (1/m, 1)$.
\end{thm}

It is not difficult to verify  (\ref{condition})  in general, as
we now examine some special cases.

\begin{rem}  The geometric distribution, i.e.
$p_k = a^k (1-a)$ for $0< a < 1$.  \rm Then
$$f(s) = \frac {1-a}{1- a s},\quad
\quad \mbox{ and } \quad \frac {(1-s)f' }{1-f } = \frac {1-a}{1- a
s} $$ which is clearly increasing in $s$.  With a little extra
work one can show that the {\it negative binomial distribution}
also satisfies (\ref{condition}).

For the {\it Poisson distribution} with parameter $\mu$, \rm
$f(s) = e^{s\mu - \mu}$ and $f'(s)$ = $\mu f(s)$.  Then
the derivative of $(1-s) f'/(1-f)$ is $\mu \phi(s) f/(1-f)^2$
where $\phi(s) = \mu(1-s)- 1+f$.  Since $\phi(1) $=0, and
$\phi'(s) <0$ for $s<1$, so $\phi(s) \geq 0$ for $s\leq 1$.

For the {\it Binomial distribution} $B(n,p)$, 
$f = (1-p+ s p)^n$ and 
one can similarly verify that  the derivative of $(1-s) f'/(1-f)$
is non-negative.
\end{rem}

\noindent{\it Proof of Proposition 1}. The regular tree is a
Galton-Watson tree with the degenerated distribution: $p_d = 1$.
In this case $f(s) = s^d$, $f'(s)$ = $ d s^{d-1}$ and
$$\frac {(1-s)f' }{1-f }= d \  \frac {s^{d-1}-s^d }{1-s^d }
$$
is indeed increasing for $s\in (0, 1)$. The condition of Theorem 3
is satisfied and Proposition 1 is proved.  \hfill $\square$

\section{Galton-Watson Trees}

Take a Galton-Watson tree with offspring distribution $\{p_k, k=
0, 1, 2, \cdots\}$ and consider the Bernoulli bond percolation on
the tree with retaining probability $p$.  The open cluster
containing the root $o$ is again a Galton-Watson tree with
offspring distribution
$$ \bar p_l = \sum_{r=0}^\infty p_{l+r} p^l
(1-p)^r C_{l+r}^r, \quad \mbox { for } l = 0, 1, 2, \cdots.$$ In
order to have an infinite cluster, the retaining probability $p >
1/m$ where $m= \sum_k k p_k$ is the mean number of offsprings.

We are now working with two Galton-Watson processes:  $\{p_k\}$
and   $\{\bar p_k\}$. Let $m= \sum_k p_k k$ be the mean of
offsprings, and $f(s)= \sum_k p_k s^k$ the generating function,
exclusively for the original $\{p_k\}$. Let $\rho$ be the
extinction probability related to $\{\bar p_k\}$, and $\lambda = 1
- p +p \rho$.

\begin{lem}
If $f(s) \neq s$, then $\rho$ is  decreasing and differentiable in
$p$ for $p \in (1/m,
1)$, and 
\begin{equation} \label{derivative}\frac {d \rho }{d p}= - \frac {(1-\rho)f'(\lambda)} { 1- pf'(\lambda) }.
\end{equation}
\end{lem}

\noindent {\it Proof}. By definition  $\rho$ satisfies the
following equation.
\begin{equation}\label{rho}\rho = \sum_{k=0}^\infty \bar p_k \rho^k =
\sum_{l=0}^\infty p_l (1-p + p \rho)^l = f(\lambda).\end{equation}
Notice that $\rho= (\lambda - 1 + p)/ p$, we can rewrite the
equation as
 \begin{equation}\label{p}  p =  \frac {1-\lambda} { 1-
f(\lambda)}.\end{equation}
So $p$ is continuous and differentiable in $\lambda$. Except the
degenerated case that $f(s) = s$, $d p/d\lambda < 0$ by the
convexity of $f$. Therefore $\lambda$ as the inverse function is
also continuous and differentiable in $p$. Taking the derivative
of the both sides of (\ref{p}) with respect to $p$, we get
$$\frac {d \lambda} {d p} = - \frac {1-f(\lambda)} { 1- pf'(\lambda) }.$$
This together with (\ref{rho}) implies (\ref{derivative}).
\hfill $\square$

\medskip
\noindent{\it Remark}:\ We know little of $\rho$. Even for regular
trees, we only have an estimate of $\rho$.
$$\rho (p) \leq \frac {1-p} {d^2 p^2 -p}  \quad \mbox {for } \ p\in
(1/d, 1).$$ An exceptional case is the binary tree for which
$\rho$  can be expressed explicitly as $(1-p)^2/p^2$, and is used
in Example 2.

\medskip We assume  that the open  cluster ${\C}_o$ containing
the root is infinite. As a Galton-Watson tree, $\C_o$  can be
constructed as follows, see \cite{L}.  Begin with the root which
is declared to be {\it green}. Add to the root a random number of
edges according to probability distribution $P(Y = k) = \bar p_k
(1-\rho^k)/(1-\rho)$. Declare the other end vertex of newly-added
edge {\it green} with probability $1-\rho$ and {\it red} with
probability $\rho$, independent of each other.  If all the
newly-added vertices are {red}, discard the entire assignment and
reassign {\it green/red} all over again.  For each {\it green}
vertex, repeat the same procedure. For each {\it red} vertex,
attach to it independently a random number of red vertices
according to the distribution $\hat p_k  = \bar p_k \rho^{k-1}$
for $k \geq 0$. The infinite tree consisting of green vertices is
called the {\it backbone} and a connected component of red
vertices is called a {\it bush}.
The backbone is a Galton-Watson tree generated according to the
distribution:
\begin{equation}\label{pk1}
 \tilde p_k  = \sum_{r=0}^\infty \bar p_{k+r} \rho^r
(1-\rho)^{k-1} C_{k+r}^r  \quad \mbox {for} \ k \geq 1.
\end{equation}
Notice that
\begin{equation}\label{pk} \tilde{p}_k = \sum_{n=k}^\infty p_n
C_n^k \lambda^{n-k} p^k (1-\rho)^{k-1} =
\frac{f^{(k)}(\lambda)}{k!} p^k (1-\rho)^{k-1},
\end{equation}
where $f^{(k)}$ is the $k$-th derivative of $f$.

\begin{lem} Denote by $M$ the
mean size of a bush. Then $M = 1/( {1 - \hat m})$, where $\hat m$
= $(1-\lambda) f'(\lambda)/(1 -f(\lambda))$.
\end{lem}
\noindent {\it Proof}. Let $\hat m$ be the mean number of
offsprings. Then
\begin{eqnarray*}
\hat m &  =  &  \sum_{k=0}^\infty \hat p_k  k = \sum_{k=0}^\infty \bar p_k \rho^{k-1} k \\
& =  & \sum_{k=0}^\infty \sum_{l=0}^\infty p_{k+l} p^k
(1-p)^l C_{k+l}^l \rho^{k-1} k \\
& =  & \sum_{n=0}^\infty \frac {p_n} \rho \lambda^n \sum_{k+l=n}
(\frac {p\rho} \lambda)^k
(\frac {1-p}\lambda)^l C_{k+l}^l  k \\
& =  & \sum_{n=0}^\infty \frac {p_n} \rho \lambda^n (\frac{ p
\rho} \lambda) n\\
& = &  p \sum_{n=0}^\infty {p_n}\  n   \lambda^{n-1} =  p
f'(\lambda) = \frac {(1-\lambda) f'(\lambda)}{1 -f(\lambda)}.
\end{eqnarray*}
 Then $M = \sum_{k=0}^\infty m^k = 1/( {1 -\hat m})$.
 \hfill $\square$

\begin{lem}
The speed $S(p)$ of the simple random walk on the backbone is
increasing in $p$ if  $(1-s)f'(s)/(1-f(s))$ is increasing in $s$.
It is continuous for all $p\leq 1$ and differentiable for $p\in(
1/m, 1)$.
\end{lem}
\noindent {\it Proof}. Let $S(p)$ be the speed of the simple
random walk on the backbone. By (1) and (\ref{pk}),
\begin{eqnarray*}
S(p) &  =  & \sum_{k=1}^\infty \tilde p_k \frac {k-1} {k+1} = 1 -
2 \sum_k \tilde p_k \frac 1 {k+1} \\  & = &  1- 2
\sum_{k=1}^\infty \frac 1 {k+1} \sum_{n=k}^\infty p_n C_n^k
\lambda^{n-k} p^k (1-\rho)^{k-1} \\ \end{eqnarray*}
\begin{eqnarray*} & = &  1- 2 \sum_{n=1}^\infty
p_n \sum_{k=1}^n \frac 1 {k+1}
C_n^k \lambda^{n-k} p^k (1-\rho)^{k-1} \\
& = &  1- \frac 2 {p ( 1- \rho)^2}  \sum_{n=1}^\infty \frac {p_n}
{n+1}  \sum_{k=1}^n \frac {n+1} {k+1}
C_n^k \lambda^{n-k}  (p-p\rho)^{k+1} \\
& = &  1- \frac 2 {p ( 1- \rho)^2}  \sum_{n=1}^\infty \frac {p_n}
{n+1} [(\lambda + p -p \rho)^{n+1} -  \lambda^{n+1} - (n+1) \lambda^{n} (p - p \rho)] \\
& = &  1- \frac 2 {p ( 1- \rho)^2}  \sum_{n=0}^\infty \frac {p_n}
{n+1} [1 -  \lambda^{n+1} - (n+1) \lambda^{n} p (1-\rho)] \\
& = &   \frac {1+\rho} {1- \rho}
 - \frac 2 {( 1- \rho)^2 p} \sum_{n=0}^\infty p_n \frac {1
 - \lambda^{n+1}} {n+1}.
 \end{eqnarray*}
The continuity and differentiability now follow from Lemma 2.1.
Take the derivative of $S(p)$ with respect to $p$, using
(\ref{derivative}) in the third equation below.
\begin{eqnarray*}
S'(p) & = &   \frac {2\rho'} {(1- \rho)^2}
 + \frac 2 {( 1- \rho)^2 p} \sum_{n=0}^\infty p_n  \lambda^n
 \lambda '
  - 2 \frac {2 p \rho' - (1-\rho)} {( 1- \rho)^3 p^2} \sum_{n=0}^\infty p_n \frac {1
 - \lambda^{n+1}} {n+1} \\
& = &   \frac {2} {(1- \rho)^2} \left[\rho'
 + \frac \rho p ( -1+ \rho + p\rho')
  -  \frac {2 p \rho' - (1-\rho)} {( 1- \rho) p^2} \sum_{n=0}^\infty p_n \frac {1
 - \lambda^{n+1}} {n+1}\right].\\
 &  = &   \frac 2 {(1- \rho)^2} \left[-(1+\rho) \frac
{(1-\rho)f'} { 1- p f' } - \frac {\rho( 1- \rho)} p  +
   \frac {1 + p f'} {( 1- p f') p^2} \sum_{n=0}^\infty p_n \frac {1
 - \lambda^{n+1}} {n+1}\right]\\
 & =  & \frac { 2 \Psi(p)}{(1- \rho)^2( 1- p f') p^2};
  \end{eqnarray*}
where $$\Psi(p)= - (1-\rho) p (\rho + p f')+ (1+pf')
\sum_{n=0}^\infty {p_n}  \frac {1 - \lambda^{n+1}}{n+1}.$$
 Furthermore, $\Psi(1/m) = 0 $ and
$$\Psi'(p) = \left[p(1-\rho)-  \sum_n
p_n \frac{1 - \lambda^{n+1}} {n+1}\right](- f' - p f''
\lambda').$$ It is enough to show $\Psi'(p) \geq 0$ for $p\in
(1/m, 1)$. This can be done in two parts. First, 
$$p(1-\rho)- \sum_n
p_n \frac{1 - \lambda^{n+1}} {n+1}  \geq 0,$$ since its value at
$p= 1/m$ is 0 and its derivative with respect to $p$ is
$(1-\rho)^2/(1- p f'(\lambda))\geq 0$. Secondly, by
(\ref{condition}),
 $$- f' - p f'' \lambda' = - \frac d {dp} (p f'(\lambda)) =
 \frac d {d\lambda} \left(\frac
{(1-\lambda)f'(\lambda)}{1-f(\lambda)}\right) \left(- \frac {d
\lambda}{dp}\right) \geq 0.$$\hfill $\square$


\medskip
 A trajectory of the simple random walk on ${\cal C}_o$
is a sequence of red vertices and green vertices.  A consecutive
sequence of red vertices, together with the two green vertices
immediately before and after the sequence, is called an {\it
excursion}. If the end of one excursion is the start of another
excursion in a trajectory, we say two excursions are consecutive.
Let $N(p, k)$ be the conditional expectation of the number of
consecutive excursions into bushes before moving to an adjacent
vertex of $x$ in the backbone, given the degree of $x$ in the
backbone is $k+1$.

\begin{lem}
$$N(p, k) =   \frac {p
  \rho}{(k+1)} \frac {f^{(k+1)} (\lambda)} {f^{(k)} (\lambda)}.$$
  Here ${f^{(k)} (\lambda)}$ is the $k$-th derivative of $f$.
\end{lem}

\noindent {\it Proof}. Suppose that the degree of vertex $x$ in
the backbone is $k+1$ and that there are $l$ bushes attached to
$x$. In this case the mean number of consecutive excursions into
bushes before moving to an adjacent vertex of the backbone is
$l/(k+1)$.

The probability that there are $l$ bushes attached to a vertex of
degree $k+1$ in the backbone is $\bar p_{k+l} (1-\rho)^{k-1}
\rho^l C_{k+l}^l$. Here we have implicitly assumed that the vertex
is not the root, since the simple random walk in a Galton-Watson
tree is transient. Take the conditional expectation of the number
of consecutive excursions into bushes.
\begin{eqnarray*}
N(p, k) & = & \frac {\sum_{l=0}^\infty \frac l {k+1} \bar p_{k+l}
 (1-\rho)^{k-1} \rho^l C_{k+l}^l} {\sum_{l=0}^\infty  \bar p_{k+l}
 (1-\rho)^{k-1} \rho^l C_{k+l}^l}\\
 & = & \frac {\sum_{l=0}^\infty \frac l {k+1} \sum_{r=0}^\infty p_{k+l+r}
  p^{k+l}(1-p)^r C_{k+l+r}^r
 (1-\rho)^k \rho^l C_{k+l}^l} {\sum_{l=0}^\infty  \sum_{r=0}^\infty p_{k+l+r}
  p^{k+l}(1-p)^r C_{k+l+r}^r
 (1-\rho)^k \rho^l C_{k+l}^l}\\
 & = & \frac {\sum_{l=1}^\infty \sum_{r=0}^\infty \frac {(k+l+r)!} {(k+1)! r! (l-1)!} p_{k+l+r}
  p^l(1-p)^r  \rho^l } {\sum_{l=0}^\infty  \sum_{r=0}^\infty \frac {(k+l+r)!} {k! r! l!} p_{k+l+r}
  p^l(1-p)^r  \rho^l }\\
 & = & \frac {\sum_{s=1}^\infty p_{k+s}\frac {(k+s)!} {(k+1)! (s-1)!}
  p \rho (1-p +p\rho)^{s-1} } {\sum_{s=0}^\infty p_{k+s}\frac {(k+s)!} {k! s!}
  (1-p +p\rho)^s }\\
 & = & 
\frac {p  \rho}{(k+1)} \frac {f^{(k+1)} (\lambda)} {f^{(k)}
(\lambda)}.
\end{eqnarray*}
 \hfill $\square$

\smallskip \noindent {\it Proof of Theorem 3}.
 While the infinite
cluster $\C_o$ is decomposed as the {\it backbone} and {\it
bushes}, \srw\ on $\C_o$ is also decomposed as \srw\ on the
backbone  and excursions into bushes. Thus \srw\ on $\C_o$ can be
regarded as a delayed simple random walk on the backbone with
delays caused by excursions into bushes.

The average time spent in an excursion into a bush, by the ergodic
theorem of Markov chains, is 2 times the size of the bush (the
number of red vertices). Thus the average of random delays at a
vertex of degree $k+1$ is, by Lemmas 2.2 and 2.4,
$$ 2 \ M  \ N(p, k)
 = \frac 2 {1-\hat m }  \ \frac {p \rho} {(k+1)}
   \ \frac {f^{(k+1)} (\lambda)} {f^{(k)}
  (\lambda)}.
$$
According to \cite{LPP}, the frequency that the simple random walk
 on the backbone visits a vertex of degree $k+1$ is $\tilde p_k $. 
By (\ref{pk}), the average of random delays at a vertex is as
follows.
 \begin{eqnarray*}
 \sum_{k=1}^\infty \tilde{p}_k 2MN(k)& =  &\frac{2 p
\rho}{1-\hat m}   \sum_{k=1}^\infty \tilde{p}_k \frac{1}{k+1}
\frac{f^{(k+1)}}{f^{(k)}}
 \\ & = &
\frac{2 p \rho}{1-pf'} \sum_{k=1}^\infty \frac{f^{(k)}}{k!} \cdot
p^k (1-\rho)^{k-1} \frac{1}{k+1} \frac{f^{(k+1)}}{f^{(k)}}
 \\ 
  & = &
\frac{2 \rho}{1-pf'} \sum_{k=1}^\infty \frac{f^{(k+1)}}{(k+1)!}
\cdot p^{k+1} (1-\rho)^{k-1}
 \\ & = &
\frac{2  \rho}{1-pf'} \frac{1}{(1-\rho)} \sum_{k=2}^\infty
\frac{f^{(k)}}{k!} \cdot p^k (1-\rho)^{k-1}
 \\ & = &
\frac{2  \rho}{1-pf'} \frac{1}{(1-\rho)} \sum_{k=2}^\infty
\tilde{p}_k
 \\ & = &
\frac{2 \rho}{1-pf'}\ \ \frac{1- \tilde{p}_1}{1-\rho}
 \\ & = &
\frac{2  \rho}{1-pf'}\ \  \frac{1- f'p}{1-\rho}
 = \frac{2 \rho}{1-\rho}.
 \end{eqnarray*}
The speed of \srw\ on $\C_o$, by the ergodic theorem of Markov
chains, is
\begin{equation}\label{ergodic}
\left (1+ \frac {2 \rho}{1-\rho}\right)^{-1} S(p) =
 \frac{1- \rho}{1+\rho} S(p)
\end{equation}
It is monotone by Lemmas 2.1 and 2.3.  Its continuity and
differentiability follow from Lemma 2.1 and (\ref{ergodic}).
\hfill $\square$

\medskip

\small\baselineskip=0.7\baselineskip \noindent{\bf
Acknowledgement}: {\it We would like to thank Yueyun Hu for
stimulating discussions. This paper is completed while the first
author is visiting MSRI, its generous support is gratefully
acknowledged.}

\noindent LMAM, School of Mathematical Sciences, Peking
University, Beijing
100871, China\\
E-mail: dayue@math.pku.edu.cn, zhangfxi@math.pku.edu.cn\\
2005-03-23
\end{document}